\documentclass[12pt,leqno]{article}
\def\diam{\mathop{\rm diam}}

\begin{document}

\title{Some Remarks Concerning Potentials on Different Spaces}

\author{S. Semmes\thanks{A lecture based on this paper was given
at the conference ``Heat kernels and analysis on manifolds''
at the Institut Henri Poincar\'e, May, 2002.}}

\date{}

\maketitle

	Let $n$ be a positive integer greater than $1$, and consider
the potential operator $P$ acting on functions on ${\bf R}^n$ defined
by
\begin{equation}
\label{def of P on R^n}
	P(f)(x) = \int_{{\bf R}^n} \frac{1}{|x-z|^{n-1}} \, f(z) \, dz.
\end{equation}
Here $dz$ denotes Lebesgue measure on ${\bf R}^n$.  More precisely,
if $f$ lies in $L^q({\bf R}^n)$, then $P(f)$ is defined almost
everywhere on ${\bf R}^n$ if $1 \le q < n$, it is defined almost
everywhere modulo constants when $q = n$, and it is defined modulo
constants everywhere if $n < q < \infty$.  (If $q = \infty$, then
one can take it to be defined modulo affine functions.)  We shall
review the reasons behind these statements in a moment.

	The case where $n = 1$ is a bit different and special, and we
shall not pay attention to it in these notes for simplicity.
Similarly, we shall normally restrict our attention to functions in
$L^q$ with $1 < q < \infty$.

	A basic fact about this operator on ${\bf R}^n$ is that if $f
\in L^q({\bf R}^n)$, then the first derivatives of $P(f)$, taken in
the sense of distributions, all lie in $L^q({\bf R}^n)$, as long as $1
< q < \infty$.  Indeed, the first derivatives of $P(f)$ are given by
first Riesz transforms of $f$ (modulo normalizing constant factors),
and these are well-known to be bounded on $L^q$ when $1 < q < \infty$.
(In connection with these statements, see \cite{St, SW}.)

	One might rephrase this as saying that $P$ maps $L^q$ into the
Sobolev space of functions on ${\bf R}^n$ whose first derivatives 
lie in $L^q$ when $1 < q < \infty$.  Instead of taking derivatives,
one can look at the oscillations of $P(f)$ more directly, as follows.
Let $r$ be a positive real number, which represents the scale at which
we shall be working.  Consider the expression
\begin{equation}
\label{frac{P(f)(x) - P(f)(y)}{r}}
	\frac{P(f)(x) - P(f)(y)}{r}.
\end{equation}

	To analyze this, let us decompose $P(f)$ into local and
distant parts at the scale of $r$.  Specifically, define operators
$L_r$ and $J_r$ by
\begin{equation}
\label{def of L_r}
	L_r(f)(x) = \int_{\{z \in {\bf R}^n : \, |z-x| < r\}}
			\frac{1}{|x-z|^{n-1}} \, f(z) \, dz
\end{equation}
and
\begin{equation}
\label{def of J_r}
	J_r(f)(x) = \int_{\{z \in {\bf R}^n : \, |z-x| \ge r\}}
			\frac{1}{|x-z|^{n-1}} \, f(z) \, dz.
\end{equation}
Thus $P(f) = L_r(f) + J_r(f)$, at least formally (we shall say more
about this in a moment), so that
\begin{equation}
\label{frac{P(f)(x) - P(f)(y)}{r} = ..., 1}
	\frac{P(f)(x) - P(f)(y)}{r} = 
	    \frac{L_r(f)(x) - L_r(f)(y)}{r} + \frac{J_r(f)(x) - J_r(f)(y)}{r}.
\end{equation}

	More precisely, $L_r(f)(x)$ is defined almost everywhere in
$x$ when $f \in L^q({\bf R}^n)$ and $1 \le q \le n$, and it is defined
everywhere when $q > n$.  These are standard results in real analysis
(as in \cite{St}), which can be derived from Fubini's theorem and
H\"older's inequality.  On the other hand, if $1 \le q < n$, then
$J_r(f)(x)$ is defined everywhere on ${\bf R}^n$, because H\"older's
inequality can be used to show that the integral converges.  This does
not work when $q \ge n$, but in this case one can consider the
integral which formally defines the difference $J_r(f)(x) - J_r(y)$.
Namely,
\begin{eqnarray}
\label{J_r(f)(x) - J_r(f)(y), 1}
\lefteqn{J_r(f)(x) - J_r(f)(y) = } 					\\
	& & \int_{{\bf R}^n}
     \biggl(\frac{1}{|x-z|^{n-1}} \, {\bf 1}_{{\bf R}^n \backslash B(x,r)}(z)
  -   \frac{1}{|y-z|^{n-1}} \, {\bf 1}_{{\bf R}^n \backslash B(y,r)}(z) \biggr)
		\, f(z) \, dz.					\nonumber
\end{eqnarray}
Here ${\bf 1}_A(z)$ denotes the characteristic function of a set $A$,
so that it is equal to $1$ when $z \in A$ and to $0$ when $z$ is not
in $A$, and $B(x,r)$ denotes the open ball in ${\bf R}^n$ with center
$x$ and radius $r$.  The integral on the right side of (\ref{J_r(f)(x)
- J_r(f)(y), 1}) does converge when $f \in L^q({\bf R}^n)$ and $q <
\infty$, because the kernel against which $f$ is integrated is bounded
everywhere, and decays at infinity in $z$ like $O(|z|^{-n})$.  This is
easy to check.

	Using this, one gets that $J_r(f)$ is defined ``modulo
constants'' on ${\bf R}^n$ when $f \in L^q({\bf R}^n)$ and $n \le q <
\infty$.  This is also why $P(f)$ can be defined modulo constants on
${\bf R}^n$ in this case (almost everywhere when $q = n$), because of
what we know about $L_r(f)$.  Note that $J_r(f)$ for different values
of $r$ can be related by the obvious formulae, with the differences
given by convergent integrals.  Using this one can see that the
definition of $P(f)$ in terms of $J_r(f)$ and $L_r(f)$ does not depend
on $r$.

	Now let us use (\ref{frac{P(f)(x) - P(f)(y)}{r} = ..., 1})
to estimate $r^{-1} (P(f)(x) - P(f)(y))$.  Specifically, in keeping
with the idea that $P(f)$ should be in the Sobolev space corresponding
to having its first derivatives be in $L^q({\bf R}^n)$ when $f$ is
in $L^q({\bf R}^n)$, $1 < q < \infty$, one would like to see that
\begin{equation}
\label{local average of the difference quotient}
	\frac{1}{|B(x,r)|} \int_{B(x,r)} \frac{|P(f)(x) - P(f)(y)|}{r}
								\, dy
\end{equation}
lies in $L^q({\bf R}^n)$, with the $L^q$ norm bounded uniformly over
$r > 0$.  Here $|A|$ denotes the Lebesgue measure of a set $A$ in
${\bf R}^n$, in this case the ball $B(x,r)$.  In fact, one can even
try to show that the supremum over $r > 0$ of (\ref{local average of
the difference quotient}) lies in $L^q$.  By well-known results, if $q
> 1$, then both conditions follow from the information that the
gradient of $P(f)$ lies in $L^q$ on ${\bf R}^n$, and both conditions
imply that the gradient of $P(f)$ lies in $L^q$.  (Parts of this work
for $q = 1$, and there are related results for the other parts.)  We
would like to look at this more directly, however.

	For the contributions of $L_r(f)$ in (\ref{frac{P(f)(x) -
P(f)(y)}{r} = ..., 1}) to (\ref{local average of the difference
quotient}), one can obtain estimates like the ones just mentioned by
standard means.  For instance, $\sup_{r > 0} r^{-1} \, L_r(f)(x)$ can
be bounded (pointwise) by a constant times the Hardy--Littlewood
maximal function of $f$ (by analyzing it in terms of sums or integrals
of averages of $f$ over balls centered at $x$).  Compare with
\cite{St, SW}.  One also does not need the fact that one has a
difference $L_r(f)(x) - L_r(f)(y)$ in (\ref{frac{P(f)(x) - P(f)(y)}{r}
= ..., 1}), but instead the two terms can be treated independently.
The localization involved is already sufficient to work back to $f$ in
a good way.

	For the $J_r(f)$ terms one should be more careful.  In
particular, it is important that we have a difference $J_r(f)(x) -
J_r(f)(y)$, rather than trying to deal with the two terms separately.
We have seen an aspect of this before, with simply having the
difference be well-defined when $f$ lies in $L^q({\bf R}^n)$ and $n
\le q < \infty$.

	Consider the auxiliary operator $T_r(f)$ defined by
\begin{equation}
\label{def of T_r}
	T_r(f)(x) = \int_{\{z \in {\bf R}^n : \, |z-x| \ge r\}}
			\frac{x-z}{|x-z|^{n+1}} \, f(z) \, dz.
\end{equation}
This is defined everywhere on ${\bf R}^n$ when $f$ lies in $L^q({\bf
R}^n)$ and $1 \le q < \infty$, because of H\"older's inequality.  Note
that $T_r(f)$ takes values in vectors, rather than scalars, because of
the presence of $x-z$ in the numerator in the kernel of the operator.
In fact,
\begin{equation}
    \nabla_x \frac{1}{|x-z|^{n-1}} = - (n-1) \frac{x-z}{|x-z|^{n+1}}.
\end{equation}
Using this and some calculus (along the lines of Taylor's theorem),
one can get that
\begin{eqnarray}
\label{estimate for J_r(f)(x) - J_r(f)(y) - (n-1) (y-x) cdot T_r(f)(x)}
\lefteqn{r^{-1} \, |J_r(f)(x) - J_r(f)(y) - (n-1) (y-x) \cdot
  T_r(f)(x)|}  \\
    & & \qquad\qquad 
 	\le C \int_{{\bf R}^n} \frac{r}{|x-z|^{n+1} + r^{n+1}} 
					\, |f(z)| \, dz	
								\nonumber
\end{eqnarray}
for a suitable constant $C$ and all $x, y \in {\bf R}^n$ with $|x-y|
\le r$.  (In other words, the kernel on the right side of
(\ref{estimate for J_r(f)(x) - J_r(f)(y) - (n-1) (y-x) cdot
T_r(f)(x)}) corresponds to the second derivatives of the kernel of
$J_r$, while $T_r$ reflects the first derivative.)

	The contribution of the right-hand side of (\ref{estimate for
J_r(f)(x) - J_r(f)(y) - (n-1) (y-x) cdot T_r(f)(x)}) to (\ref{local
average of the difference quotient}) satisfies the kind of estimates
that we want, by standard results.  (The right-hand side of
(\ref{estimate for J_r(f)(x) - J_r(f)(y) - (n-1) (y-x) cdot
T_r(f)(x)}) is approximately the same as the Poisson integral of
$|f|$.  Compare with \cite{St, SW} again.)  The remaining piece to
consider is
\begin{equation}
	(n-1) \, r^{-1} \, (y-x) \cdot T_r(f)(x).
\end{equation}
After averaging in $y$ over $B(x,r)$, as in (\ref{local average of the
difference quotient}), we are reduced to looking simply at
$|T_r(f)(x)|$.  Here again the Riesz transforms arise, but in the form
of the truncated singular integral operators, rather than the singular
integral operators themselves (with the limit as $r \to 0$).  By
well-known results, these truncated operators $T_r$ have the property
that they are bounded on $L^q({\bf R}^n)$ when $1 < q < \infty$, with
the operator norm being uniformly bounded in $r$.  Moreover, the
maximal truncated operator
\begin{equation}
	\sup_{r > 0} |T_r(f)(x)|
\end{equation}
is bounded on $L^q({\bf R}^n)$, $1 < q < \infty$.  See \cite{St, SW}.

	These statements are all closely related to the original one
concerning the way that the first derivatives of $P(f)$ are given by
first Riesz transforms of $f$ (up to constant multiples), and lie in
$L^q({\bf R}^n)$ when $f$ does and $1 < q < \infty$.  Instead of
comparing the derivatives of $P(f)$ with Riesz transforms of $f$, we
compare oscillations of $P(f)$ at the scale of $r$ with averages of
$f$ and truncated Riesz transforms of $f$ at the scale of $r$.  We do
this directly, rather than going through derivatives and integrations
of them.

	A nice feature of this discussion is that it lends itself in a
simple manner to more general settings.  In particular, it applies to
situations in which it may not be as convenient to work with
derivatives and integrations of them, while measurements of
oscillations at the scale of $r$ and related estimates still make
sense.

	Instead of ${\bf R}^n$, let us consider a set $E$ in some
${\bf R}^m$.  Let us assume that $E$ is \emph{Ahlfors-regular of
dimension $n$}, by which we mean that $E$ is closed, has at least two
elements (to avoid degeneracies), and that there is a constant $C > 0$
such that
\begin{equation}
\label{Ahlfors-regularity condition}
	C^{-1} \, t^n \le H^n(E \cap \overline{B}(x,t)) \le C \, t^n
\end{equation}
for all $x \in E$ and $t > 0$ with $t \le \diam E$.  Here $H^n$
denotes $n$-dimensional Hausdorff measure (as in \cite{Fe, Ma}), and
$\overline{B}(x,t)$ denotes the closed ball in the ambient space ${\bf
R}^m$ with center $x$ and radius $t$.

	This condition on $E$ ensures that $E$ behaves
measure-theoretically like ${\bf R}^n$, even if it could be very
different geometrically.  Note that one can have Ahlfors-regular sets
of noninteger dimension, and in fact of any dimension in $(0,m]$ (for
subsets of ${\bf R}^m$).

	Given a function $f$ on $E$, define $P(f)$ on $E$ in the same
manner as before, i.e., by
\begin{equation}
\label{def of P on E}
	P(f)(x) = \int_{E} \frac{1}{|x-z|^{n-1}} \, f(z) \, dz,
\end{equation}
where now $dz$ denotes the restriction of $H^n$-measure to $E$.  Also,
$|x-z|$ uses the ordinary Euclidean distance on ${\bf R}^m$.

	The Ahlfors-regularity of dimension $n$ of $E$ ensures that
$P(f)$ has many of the same basic properties on $E$ as on ${\bf R}^n$.
In particular, if $f$ is in $L^q(E)$, then $P(f)$ is defined almost
everywhere on $E$ (using the measure $H^n$ still) when $1 \le q < n$,
it is defined almost everywhere modulo constants on $E$ when $q = n$,
and it is defined everywhere on $E$ modulo constants when $n < q <
\infty$.  One can show these statements in essentially the same manner
as on ${\bf R}^n$, and related results about integrability, bounded
mean oscillation, and H\"older continuity can also be proven in
essentially the same manner as on ${\bf R}^n$.  

	What about the kind of properties discussed before, connected
to Sobolev spaces?  For this again one encounters operators on
functions on $E$ with kernels of the form
\begin{equation}
\label{kernels for operators}
	\frac{x-z}{|x-z|^{n+1}}.
\end{equation}
It is not true that operators like these have the same kind of
$L^q$-boundedness properties as the Riesz transforms do for arbitrary
Ahlfors-regular sets in ${\bf R}^m$, but this is true for integer
dimensions $n$ and ``uniformly rectifiable'' sets $E$.  In this
connection, see \cite{Ca, CDM, CMM, Da1, Da2, Da3, Da4, DS1, DS2, Ma,
MMV}, for instance (and further references therein).

	When $E$ is not a plane, the operators related to the kernels
(\ref{kernels for operators}) are no longer convolution operators, and
one loses some of the special structure connected to that.  However,
many real-variable methods still apply, or can be made to work.  See
\cite{CW1, CW2, CM, Jo}.  For example, the Hardy--Littlewood maximal
operator still behaves in essentially the same manner as on Euclidean
spaces, as do various averaging operators (as were used in the earlier
discussion).  Although one does not know that singular integral
operators with kernels as in (\ref{kernels for operators}) are bounded
on $L^q$ spaces for arbitrary Ahlfors-regular sets $E$, there are
results which say that boundedness on one $L^q$ space implies
boundedness on all others, $1 < q < \infty$.  Boundedness of singular
integral operators (of the general Calder\'on--Zygmund type) implies
uniform boundedness of the corresponding truncated integral operators,
and also boundedness of the maximal truncated integral operators.

	At any rate, a basic statement now is the following.  Let $n$
be a positive integer, and suppose that $E$ is an Ahlfors-regular set
in some ${\bf R}^m$ which is ``uniformly rectifiable''.  Define the
potential operator $P$ on functions on $E$ as in (\ref{def of P on
E}).  Then $P$ takes functions in $L^q(E)$, $1 < q < \infty$, to
functions on $E$ (perhaps modulo constants) which satisfy ``Sobolev
space'' conditions like the ones on ${\bf R}^n$ for functions with
gradient in $L^q$.  In particular, one can look at this in terms of
$L^q$ estimates for the analogue of (\ref{local average of the
difference quotient}) on $E$, just as before.  These estimates can be
derived from the same kinds of computations as before, with averaging
operators and operators like $T_r$ in (\ref{def of T_r}), but now on
$E$.  The estimates for $T_r$ use the assumption of uniform
rectifiability of $E$ (boundedness of singular integral operators).
The various other integral operators, with the absolute values inside
the integral sign, are handled using only the Ahlfors-regularity of
$E$.

	Note that for sets $E$ of this type, one does not necessarily
have the same kind of properties concerning integrating derivatives as
on ${\bf R}^n$.  In other words, one does not automatically get as
much from looking at infinitesimal oscillations, along the lines of
derivatives, as one would on ${\bf R}^n$.  The set $E$ could be quite
disconnected, for instance.  However, one gets the same kind of
estimates at larger scales for the potentials that one would normally
have on ${\bf R}^n$ for a function with its first derivatives in
$L^q$, by looking at a given scale $r$ directly (rather than trying to
integrate bounds for infinitesimal oscillations), as above.

	For some topics related to Sobolev-type classes on general
spaces, see \cite{FHK, Ha, HaK1, HaK2, HeK1, HeK2} (and references
therein).

	Although the potential operator in (\ref{def of P on E}) has a
nice form, it is also more complicated than necessary.  Suppose that
$E$ is an $n$-dimensional Lipschitz graph, or that $E$ is simply
bilipschitz--equivalent to ${\bf R}^n$, or to a subset of ${\bf R}^n$.
In these cases the basic subtleties for singular integral operator
with kernel as in (\ref{kernels for operators}) already occur.
However, one can obtain potential operators with the same kind of nice
properties by making a bilipschitz change of variables into ${\bf
R}^n$, and using the classical potential operator there.  This leads
back to the classical first Riesz transforms on ${\bf R}^n$, as in
\cite{St, SW}.

	Now let us consider a rather different kind of situation.
Suppose that $E$ is an Ahlfors-regular subset of dimension $n$ of some
${\bf R}^m$ again.  For this there will be no need to have particular
attention to integer values of $n$.  Let us say that $E$ is a
\emph{snowflake} of order $\alpha$, $0 < \alpha < 1$, if there is a
constant $C_1$ and a metric $\rho(x,y)$ on $E$ such that
\begin{equation}
\label{snowflake condition}
	C_1^{-1} \, |x-y| \le \rho(x,y)^\alpha \le C_1 \, |x-y|
\end{equation}
for all $x, y \in E$.

	In this case, let us define a potential operator
$\widetilde{P}$ on functions on $E$ by
\begin{equation}
\label{def of widetilde{P}}
	\widetilde{P}(f)(x) 
	  = \int_E \frac{1}{\rho(x,z)^{\alpha (n-1)}} \, f(z) \, dz.
\end{equation}
Here $dz$ denotes the restriction of $n$-dimensional Hausdorff measure
to $E$ again.  This operator is very similar to the one before, since
$\rho(x,z)^{\alpha (n-1)}$ is bounded from above and below by constant
multiples of $|x-z|^{n-1}$, so that the kernel of $\widetilde{P}$ is
bounded from above and below by constant multiples of the kernel of
the operator $P$ in (\ref{def of P on E}).

	This operator enjoys the same basic properties as before, with
$\widetilde{P}(f)$ being defined almost everywhere when $f$ lies in
$L^q(E)$ and $1 \le q < n$, defined modulo constants almost everywhere
when $q = n$, and defined modulo constants everywhere when $n < q <
\infty$, for essentially the same reasons as in the previous
circumstances.  However, there is a significant difference with this
operator, which one can see as follows.  Let $x$, $y$, $z$ be three
points in $E$, with $x \ne z$ and $y \ne z$.  Then
\begin{equation}
\label{inequality for differences of kernels}
 	\biggl| \frac{1}{\rho(x,z)^{\alpha (n-1)}}  
                  - \frac{1}{\rho(y,z)^{\alpha (n-1)}} \biggr|
 \le C \, \frac{\rho(x,y)}{\min(\rho(x,z),\rho(y,z))^{\alpha (n-1) + 1}}
\end{equation}
for some constant $C$ which does not depend on $x$, $y$, or $z$,
but only on $\alpha (n-1)$.  Indeed, one can choose $C$ so that
\begin{equation}
\label{inequality simply for positive numbers}
	\bigl| a^{\alpha (n-1)} - b^{\alpha (n-1)} \bigr|
		\le C \, \frac{|a-b|}{\min(a,b)^{\alpha (n-1) + 1}}
\end{equation}
whenever $a$ and $b$ are positive real numbers.  This is an elementary
observation, and in fact one can take $C = \alpha (n-1)$.  One can get
(\ref{inequality for differences of kernels}) from (\ref{inequality
simply for positive numbers}) by taking $a = \rho(x,z)$ and $b =
\rho(y,z)$, and using the fact that
\begin{equation}
	| \rho(x,z) - \rho(y,z)| \le \rho(x,y).
\end{equation}
This last comes from the triangle inequality for $\rho(\cdot, \cdot)$,
which we assumed to be a metric.

	Using the snowflake condition (\ref{snowflake condition}), we
can obtain from (\ref{inequality for differences of kernels}) that
\begin{equation}
\label{inequality for differences of kernels, 2}
 	\biggl| \frac{1}{\rho(x,z)^{\alpha (n-1)}}  
                  - \frac{1}{\rho(y,z)^{\alpha (n-1)}} \biggr|
  \le C' \, \frac{|x-y|^{1/\alpha}}{\min(|x-z|, |y-z|)^{(n-1) + 1/\alpha}}
\end{equation}
for all $x, y, z \in {\bf R}^n$ with $x \ne z$, $y \ne z$, and with a
modestly different constant $C'$.  The main point here is that the
exponent in the denominator on the right side of the inequality is
strictly larger than $n$, because $\alpha$ is required to lie in
$(0,1)$.  In the previous contexts, using the kernel $1/ |x-z|^{n-1}$
for the potential operator, there was an analogous inequality with
$\alpha = 1$, so that the exponent in the denominator was equal to
$n$.

	With an exponent larger than $n$, there is no need for
anything like singular integral operators here.  More precisely, there
is no need for the operators $T_r$ in (\ref{def of T_r}) here; one can
simply drop them, and estimate the analogue of $|J_r(f)(x) -
J_r(f)(y)|$ when $|x - y| \le r$ directly, using (\ref{inequality for
differences of kernels, 2}).  In other words, one automatically gets
an estimate like (\ref{estimate for J_r(f)(x) - J_r(f)(y) - (n-1)
(y-x) cdot T_r(f)(x)}) in this setting, without the $T_r$ term, and
with some minor adjustments to the right-hand side.  Specifically, the
$r$ in the numerator on the right side of (\ref{estimate for J_r(f)(x)
- J_r(f)(y) - (n-1) (y-x) cdot T_r(f)(x)}) would become an
$r^{1/\alpha -1}$ in the present situation, and the exponent $n+1$ in
the denominator would be replaced with $n - 1 + 1/\alpha$.  This leads
to the same kinds of results in terms of $L^q$ norms and the like as
before, because the rate of decay is enough so that the quantities in
question still look like suitable averaging operators in $f$.  (That
is, they are like Poisson integrals, but with somewhat less decay.
The decay is better than $1/|x-z|^n$, which is the key.  As usual, see
\cite{St, SW} for similar matters.)

	The bottom line is that if we use the potential operator
$\widetilde{P}$ from (\ref{def of widetilde{P}}) instead of the
operator $P$ from (\ref{def of P on E}), then the two operators are
approximately the same in some respects, with the kernels being of
comparable size in particular, but in this situation the operator
$\widetilde{P}$ has the nice feature that it automatically enjoys the
same kind of properties as in the ${\bf R}^n$ case, in terms of
estimates for expressions like (\ref{local average of the difference
quotient}) (under the snowflake assumption for $E$).  That is, one
automatically has that $\widetilde{P}(f)$ behaves like a function in a
Sobolev class corresponding to first derivatives being in $L^q$ when
$f$ lies in $L^q$.  One does not need $L^q$ estimates for singular
integral operators for this, as would arise if we did try to use the
operator $P(f)$ from (\ref{def of P on E}).

	These remarks suggest numerous questions...

	Of course, some other basic examples involve nilpotent Lie
groups, like the Heisenberg group, and their invariant geometries.

	As a last comment, note that for the case of snowflakes we
never really needed to assume that $E$ was a subset of some ${\bf
R}^m$.  One could have worked just as well with abstract metric spaces
(still with the snowflake condition).  However, Assouad's embedding
theorem \cite{A1, A2, A3} provides a way to go back into some ${\bf
R}^m$ anyway.  The notion of uniform rectifiability makes sense for
abstract metric spaces, and not just subsets of ${\bf R}^m$, and an
embedding into some ${\bf R}^m$ is sometimes convenient.  In this
regard, see \cite{Se}.

\end{document}